\documentstyle[twoside]{article}

\oddsidemargin=\evensidemargin
\catcode`\@=11
\long\def\@makefntext#1{
\protect\noindent \hbox to 3.2pt {\hskip-.9pt  
$^{{\eightrm\@thefnmark}}$\hfil}#1\hfill}		

\def\ps@myheadings{\let\@mkboth\@gobbletwo		
\def\@oddhead{\hbox{}
\rightmark\hfil\eightrm\thepage}   
\def\@oddfoot{}\def\@evenhead{\eightrm\thepage\hfil
\leftmark\hbox{}}\def\@evenfoot{}
\def\sectionmark##1{}\def\subsectionmark##1{}}

\catcode`@=11		      		     
\def\ps@plain{\let\@mkboth\@gobbletwo
     \def\@oddhead{}\def\@oddfoot{\eightrm\hfil\thepage
     \hfil}\def\@evenhead{}\let\@evenfoot\@oddfoot}

\renewcommand{\thefootnote}{\fnsymbol{footnote}}

\newcounter{sectionc}\newcounter{subsectionc}\newcounter{subsubsectionc}
\renewcommand{\section}[1] {\vspace{12pt}\addtocounter{sectionc}{1} 
\setcounter{subsectionc}{0}\setcounter{subsubsectionc}{0}\noindent 
	{\tenbf\thesectionc. #1}\par\vspace{5pt}}
\renewcommand{\subsection}[1] {\vspace{12pt}\addtocounter{subsectionc}{1} 
	\setcounter{subsubsectionc}{0}\noindent 
	{\bf\thesectionc.\thesubsectionc. 
	{\kern1pt \bfit #1}}\par\vspace{5pt}}
\renewcommand{\subsubsection}[1] {\vspace{12pt}
	\addtocounter{subsubsectionc}{1}
	\noindent
	{\tenrm\thesectionc.\thesubsectionc.\thesubsubsectionc.	{\kern1pt 
	\it #1}}\par\vspace{5pt}}
\newcommand{\nonumsection}[1] {\vspace{12pt}\noindent{\tenbf #1}
	\par\vspace{5pt}}

\newcounter{appendixc}
\newcounter{subappendixc}[appendixc]
\newcounter{subsubappendixc}[subappendixc]

\renewcommand{\appendix}[1] {\vspace{12pt}	
	\refstepcounter{appendixc}		
	\setcounter{figure}{0}
	\setcounter{table}{0}
	\setcounter{lemma}{0}
	\setcounter{theorem}{0}
	\setcounter{corollary}{0}
	\setcounter{definition}{0}
	\setcounter{equation}{0}
	\renewcommand{\thefigure}{\Alph{appendixc}.\arabic{figure}}
	\renewcommand{\thetable}{\Alph{appendixc}.\arabic{table}}
	\renewcommand{\theappendixc}{\Alph{appendixc}}
	\renewcommand{\thelemma}{\Alph{appendixc}.\arabic{lemma}}
	\renewcommand{\thetheorem}{\Alph{appendixc}.\arabic{theorem}}
	\renewcommand{\thedefinition}{\Alph{appendixc}.\arabic{definition}}
	\renewcommand{\thecorollary}{\Alph{appendixc}.\arabic{corollary}}
	\renewcommand{\theequation}{\Alph{appendixc}.\arabic{equation}}
	\noindent{\tenbf Appendix \theappendixc #1}\par\vspace{5pt}}

\newcommand{\textlineskip}{\baselineskip=13pt}
\newcommand{\smalllineskip}{\baselineskip=10pt}

\newcommand{\copyrightheading}[1]
	{\vspace*{-2.5cm}\smalllineskip{\flushleft
	{\footnotesize \quad }\\
   	{\footnotesize \quad }\\
         }}

\def\abstracts#1{{
	\centering{\begin{minipage}{4.5in}\footnotesize\baselineskip=10pt
	\centerline{ABSTRACT} 
	\parindent=15pt #1\par 
	\end{minipage}}\par}} 

\def\keywords#1{{ 
	\centering{\begin{minipage}{4.5in}\footnotesize\baselineskip=10pt
	{\footnotesize\it Keywords}\/: #1
	\end{minipage}}\par}}

\renewenvironment{thebibliography}[1]
	{\frenchspacing
	 \ninerm\baselineskip=11pt
	 \begin{list}{[\arabic{enumi}]}
	{\usecounter{enumi}\setlength{\parsep}{0pt}
	 \setlength{\leftmargin 13.7pt}{\rightmargin 0pt} 
	 \setlength{\itemsep}{0pt} \settowidth
	{\labelwidth}{[#1]}\sloppy}}{\end{list}}

\newcounter{itemlistc}
\newcounter{romanlistc}
\newcounter{alphlistc}
\newcounter{arabiclistc}

\newcommand{\fcaption}[1]{
        \refstepcounter{figure}
        \setbox\@tempboxa = \hbox{\footnotesize Fig.~\thefigure. #1}
        \ifdim \wd\@tempboxa > 5in
           {\begin{center}
        \parbox{5in}{\footnotesize\smalllineskip Fig.~\thefigure. #1}
            \end{center}}
        \else
             {\begin{center}
             {\footnotesize Fig.~\thefigure. #1}
              \end{center}}
        \fi}

\newcommand{\tcaption}[1]{
        \refstepcounter{table}
        \setbox\@tempboxa = \hbox{\footnotesize Table~\thetable. #1}
        \ifdim \wd\@tempboxa > 5in
           {\begin{center}
        \parbox{5in}{\footnotesize\smalllineskip Table~\thetable. #1}
            \end{center}}
        \else
             {\begin{center}
             {\footnotesize Table~\thetable. #1}
              \end{center}}
        \fi}

\def\pmb#1{\setbox0=\hbox{#1}
	\kern-.025em\copy0\kern-\wd0
	\kern.05em\copy0\kern-\wd0
	\kern-.025em\raise.0433em\box0}

\def\fnt#1#2{\footnotetext{\kern-.3em
	{$^{\mbox{\scriptsize #1}}$}{#2}}}

\def\fpage#1{\begingroup
\voffset=.3in
\thispagestyle{empty}\begin{table}[b]\centerline{\footnotesize #1}
	\end{table}\endgroup}

\font\tenrm=cmr10
 
\font\tenbf=cmbx10
\font\bfit=cmbxti10 at 10pt
\font\ninerm=cmr9
\font\nineit=cmti9
\font\ninebf=cmbx9
\font\eightrm=cmr8

\newtheorem{thm}{Theorem}[sectionc]   

\newtheorem{lemma}{Lemma}[sectionc]
\newtheorem{corollary}{Corollary}[sectionc]
\newtheorem{definition}{Definition}[sectionc]
\newtheorem{remark}{Remark}[sectionc]
\newtheorem{example}{Example}[sectionc]  
\newtheorem{proposition}{Proposition}[sectionc]

\def\@begintheorem#1#2{\trivlist	
	\item[\hskip\labelsep{\bf #1\ #2.}]} 
\def\@opargbegintheorem#1#2#3{\trivlist
	\item[\hskip\labelsep{\bf #1\ #2\ (#3).}]}

    	{\setcounter{itemlistc}{0}		
	 \begin{list}{$\bullet$}		
	{\usecounter{itemlistc}			
	 \leftmargin10pt	       
	 \setlength{\parsep}{0pt}
	 \setlength{\itemsep}{0pt}     
	}}{\end{list}}

	{\setcounter{romanlistc}{0}		
	 \begin{list}{$($\roman{romanlistc}$)$}	
	{\usecounter{romanlistc}		
	 \leftmargin18pt 
	 \setlength{\parsep}{0pt}
	 \setlength{\itemsep}{0pt}	
	 \settowidth{\labelwidth}{#1}                          
	}}{\end{list}}

	{\setcounter{enumii}{0}			
	 \begin{list}{$($\alph{enumii}$)$}	
	{\usecounter{enumii}			
	 \leftmargin18pt		
	 \setlength{\parsep}{0pt}
	 \setlength{\itemsep}{0pt}	
	 \settowidth{\labelwidth}{#1}                          
	}}{\end{list}}

\def\qed{\hbox{${\vcenter{\vbox{			
   \hrule height 0.4pt\hbox{\vrule width 0.4pt height 6pt
   \kern5pt\vrule width 0.4pt}\hrule height 0.4pt}}}$}}

\renewcommand{\thefootnote}{\fnsymbol{footnote}}  

\def\theequation{\thesectionc.\arabic{equation}}  
\begin{document}
\setlength{\textheight}{7.7truein}  


\normalsize\textlineskip
\thispagestyle{empty}
\setcounter{page}{1}

\copyrightheading{}		    

\vspace*{0.88truein}

\fpage{1}
\centerline{\bf APPROXIMATING JONES COEFFICIENTS AND OTHER}
\baselineskip=13pt
\centerline{\bf LINK INVARIANTS BY VASSILIEV INVARIANTS}
\vspace*{0.37truein}
\centerline{\footnotesize ILYA KOFMAN}
\baselineskip=12pt
\centerline{\footnotesize\it Department of Mathematics, University of Maryland,}
\baselineskip=10pt
\centerline{\footnotesize\it College Park, MD 20742, U.S.A.}

\vspace*{10pt}
\centerline{\footnotesize YONGWU RONG}
\baselineskip=12pt
\centerline{\footnotesize\it Department of Mathematics, George Washington University}
\baselineskip=10pt
\centerline{\footnotesize\it Washington, DC 20052, U.S.A.}

\vspace*{0.225truein}
\textlineskip			
\vspace*{12pt}			

\vspace*{0.21truein} 
\abstracts{We find approximations by Vassiliev invariants for the
coefficients of the Jones polynomial and all specializations of the
HOMFLY and Kauffman polynomials.  Consequently, we obtain approximations
of some other link invariants arising from the homology of branched covers 
of links.}{}{}{}

\vspace*{10pt}

\setcounter{footnote}{0}
\renewcommand{\thefootnote}{\alph{footnote}}


\vspace*{1pt}\textlineskip	
\section{Introduction}	
\vspace*{-0.5pt}

A well-known conjecture in the theory of Vassiliev invariants is that
these invariants are dense in the space of all numerical knot invariants.
This was posed as a problem in \cite{BL} as follows:  given any numerical
knot invariant $\phi : L \to {\bf Q},$ does there exist a sequence of
Vassiliev invariants $\{v_i: L \to {\bf Q}, i=2,3,4,\ldots\}$ such that

\[ \lim_{i\to\infty} v_i(L) = \phi(L) \]

In this note, we find approximations by Vassiliev invariants for the
coefficients of the Jones polynomial and all specializations of the
HOMFLY and Kauffman polynomials.  Consequently, we obtain approximations
of some other link invariants.  This note is organized as follows:  In
Section 2, we show that every Jones coefficient is the limit of a sequence
of Vassiliev invariants.  In Section 3, given any $d$, for any Jones
polynomial of degree bounded by $d$, we find an explicit finite formula
for its coefficients in terms of Vassiliev invariants.  In Section 4, we
find an explicit infinite approximation for any Jones coefficient.  In
Section 5, we extend the results to any specialization of the HOMFLY and
Kauffman polynomials:  we find the finite formula for polynomials of
bounded degree and the infinite formula for all polynomials.  In Section
6, we find approximations by Vassiliev invariants for some link invariants
arising from the homology of branched covers of links.  In Section 7, we
discuss some conjectures related to approximations by Vassiliev
invariants.

\section{Approximating Jones coefficients by Vassiliev invariants:
Existence theorem}

Let $J_L(t)$ denote the Jones polynomial of a knot $L$.  Suppose
$ J_L(t) = a_{-m}t^{-m} + \cdots + a_0 + \cdots + a_nt^n$,
 where $a_{-m}$ and $a_n$ are nonzero.
We call the {\it degree} of the Laurent polynomial $d=\max(m,n)$.  
In \cite{BL}, it was shown that if we let $t=e^x$,

\begin{equation} \label{exp}
 J_L(e^x) = \sum_{i=0}^{\infty} \left( \frac{1}{i!}
\sum_{k=-m}^{n}k^ia_k(L) \right)x^i =   \sum_{i=0}^{\infty} v_i(L)x^i 
\end{equation}
then $v_i$ is a Vassiliev invariant of order $i$.  Henceforth, we will
refer to $\{v_i\}$ as the Vassiliev invariants obtained from the
coefficients of the expansion above.

This can be reformulated in terms of the following infinite matrix:
\begin{equation} \label{Van}
  	\left( \matrix{
	\cdots &1&1&1&1&1&1&1&\cdots\cr 
	\cdots&-3&-2&-1&0&1&2&3&\cdots\cr
	\cdots&(-3)^2&(-2)^2&(-1)^2&0&1^2&2^2&3^2&\cdots\cr
	&\vdots&\vdots&\vdots&\vdots&\vdots&\vdots&\vdots\cr
	} \right) 
	\left( \matrix{
	\vdots\cr a_{-1}\cr a_0\cr a_1 \cr  \vdots \cr 
	} \right)
=
	\left( \matrix{
	v_0 \cr v_1 \cr 2!v_2 \cr 3!v_3 \cr  \vdots \cr 
	} \right)	
\end{equation}

Recall that a Vandermonde matrix has the following form (see, e.g., \cite{F}):
\[  V(x_1,\dots,x_r) = 
	\left( \matrix{
	1&\cdots&1 \cr  
	x_1&\cdots&x_r \cr 
	x_1^2&\cdots&x_r^2 \cr 
	\vdots&&\vdots \cr 
	x_1^{r-1}&\cdots&x_r^{r-1} \cr 
	} \right) \]
\[ \det V(x_1,\dots,x_r) = \prod_{i<j}(x_j-x_i) \]
Thus, $V(x_1,\dots,x_r)$ is invertible if $x_i\neq x_j$ for all $i\neq j$.
The matrix in $(\ref{Van})$ is a Vandermonde matrix for every finite square block which
contains the first row.  From the resulting system of linear equations, we
obtain the following existence theorem:
\textlineskip
\vspace*{12pt}

\begin{thm} \label{R}
Given any knot $L$, let $J_L(t)$ be the Jones polynomial, and $a_i$ 
be its $i^{\rm th}$ coefficient.  Then for each $i$, $a_i(L)$ is the
limit of a sequence of Vassiliev invariants.
\end{thm}
\textlineskip
\vspace*{12pt}

{\bf Proof.} 
For any coefficient $a_i$, we will define a sequence of Vassiliev
invariants 
$\alpha_{1,i}, \alpha_{2,i}, \ldots$ and show that $\lim_{n\rightarrow \infty}
\alpha_{n,i} (L) = a_i(L)$.

We now let $t=e^x$ and consider the expansion (\ref{exp}).
If $n \geq d$, we obtain the following system of linear equations:

\begin{equation} \label{lin}
\left\{ \begin{array}{rclll}
a_{-n} & +  \ldots + a_{-1} + a_0 + a_1 + \ldots + & a_n & = & v_0  \\
(-n)a_{-n} & + \ldots + (-1) a_{-1} +  a_1 + \ldots  + & na_n & = & v_1  \\ 
\vdots & & & & \vdots \\
(-n)^ka_{-n} & + \ldots + (-1)^ka_{-1} +  a_1 + \ldots + & n^k a_n & = &
v_k
\end{array} 
\right.
\end{equation}

\noindent
with $2n+1$ variables $a_{-n}, \ldots, a_n$, and $k+1$ equations.
When $k=2n$, the system of linear equations has a unique solution since
the coefficient matrix is an invertible finite block of the Vandermonde
matrix from (\ref{Van}).
Denote the solution by the vector
$(\alpha_{n, -n}, \ldots, \alpha_{n, -1},  \alpha_{n, 0}, \alpha_{n, 1}, \ldots, \alpha_{n,
n})$.

We claim $\lim_{n\rightarrow \infty} \alpha_{n, i}(L)=a_i(L)$.  The claim
follows immediately from the following lemma:
\textlineskip
\vspace*{12pt}

\begin{lemma} \label{L}
For any knot $L$, $\alpha_{n,i}(L)=a_i(L)$ for all $n\geq d$.
\end{lemma}
\textlineskip
\vspace*{12pt}

\noindent
When $n\geq d$, $J_{L, n}=J_L$.  We let $\alpha_n(L)=(\alpha_{n, -n}(L), \ldots,
\alpha_{n, 0}(L),  \ldots, \alpha_{n, n}(L))$ be given by the coefficients of
$J_L(t)$.  Therefore, $\alpha_n(L)$ is the unique solution satisfying the
above system of linear equations.  If we now fix $i$, for all $n\geq d$,
it follows that $\alpha_{n,i}(L)=a_i(L)$.  This completes the proof of the
theorem. \hfill \qed\kern0.8pt

\textlineskip
\vspace*{12pt}

\begin{corollary} \label{z}
For any knot $L$ and any fixed complex number $z$, $J_L(z)$ is a limit of
Vassiliev invariants.
\end{corollary}
\textlineskip
\vspace*{12pt}

{\bf Proof.}
For each $n$, let 
\[ \alpha_n(L)=(\ldots, 0, \alpha_{n, -n}(L), \ldots, \alpha_{n, -1}(L), \alpha_{n, 0}(L),
  \alpha_{n, 1}(L), \ldots, \alpha_{n, n}(L), 0, \ldots) \] 

\noindent
be the sequence of infinite vectors as defined in Theorem \ref{R}.  Let 
\[ g_n(L)=\alpha_{n, -n}(L) z^{-n} + \ldots + \alpha_{n,0}(L) +\ldots +
\alpha_{n,n}(L)z^n. \]

\noindent
Then $g_n$ is a Vassiliev invariant since it is a linear combination of
such invariants.  By Lemma \ref{L}, when $n\geq d$, we have
$\alpha_{n, i}(L)=a_i(L)$ for all $i$.  Thus $g_n(L)=J_L(z)$.

\hfill \qed\kern0.8pt

\section{Approximating Jones coefficients by Vassiliev invariants:
Bounded degree case}

For any given $d$ (in particular, for any given knot), we can obtain
explicit solutions to the linear system (\ref{lin}) and obtain a formula
for all $\alpha_{n,i}$.  We will use (\ref{Van}) and compute the inverse of
the $(2d+1) \times (2d+1)$ Vandermonde matrix which is symmetric about the
column with zeros.  For this, we need the following generating function:

\textlineskip
\vspace*{12pt}

\begin{definition}
\[ f_{d,n}(v) = \prod_{j\neq n \atop j=-d}^{d}  \frac{v-j}{n-j} =
\frac{(-1)^{n+d}}{(d+n)! \, (d-n)!} \prod_{{j\neq n}\atop{j=-d}}^d
(v-j) \]
\end{definition}

The proof of the following proposition is immediate from the definition:

\begin{proposition} For any $m \in {\bf Z}$ such that $-d \leq m \leq d$,
\[ f_{d,n}(m) = \cases{ 1 & if $ m=n$ \cr 0 & if  $m \neq n$ \cr }\]
\end{proposition}

\begin{thm} \label{Jb}
For any Jones polynomial of a knot of degree $\leq d$, 
\[ a_n = \sum_{i=0}^{2d} f_{d,n}^{(i)}(0) v_i, \quad {\rm where }\;
f_{d,n}(v) = \prod_{j\neq n \atop j=-d}^{d}  \frac{v-j}{n-j} \]
In other words, $\alpha_{d,n} = \sum\limits_{i=0}^{2d} f_{d,n}^{(i)}(0) v_i$.
\end{thm}
\textlineskip
\vspace*{12pt}

{\bf Proof.} \qquad  Let $c_{n,j} = \frac{1}{j!} f_{d,n}^{(j)}(0)$, the $j^{th}$
coefficient of the polynomial $f_{d,n}(v)$.
\[  	\left( \matrix{
	c_{-d,0}&\cdots&c_{-d,2d} \cr  
	\vdots&&\vdots \cr 
	c_{d,0}&\cdots&c_{d,2d} \cr 
	} \right) 
	\left( \matrix{
	1&\cdots&1&\cdots&1 \cr  
	-d&\cdots&0&\cdots&d \cr 
	\vdots&&\vdots&&\vdots \cr 
	(-d)^{2d}&\cdots&0&\cdots&d^{2d} \cr 
	} \right)
= \]
\[
	\left( \matrix{
	f_{d,-d}(-d)&\cdots&f_{d,-d}(d) \cr 
	\vdots&&\vdots \cr 
	f_{d,d}(-d)&\cdots&f_{d,d}(d) \cr 
	} \right)	
= I_{(2d+1) \times (2d+1)} \]

Therefore,
\[  	\left( \matrix{
	a_{-d} \cr  \vdots \cr a_0 \cr  \vdots \cr a_d \cr 
	} \right)
=
	\left( \matrix{
	c_{-d,0}&\cdots&c_{-d,2d} \cr  
	\vdots&&\vdots \cr 
	c_{d,0}&\cdots&c_{d,2d} \cr 
	} \right)
	\left( \matrix{
	v_0 \cr v_1 \cr 2!v_2 \cr  \vdots \cr (2d)!v_{2d} \cr 
	} \right) \]
 \hfill \qed\kern0.8pt

\textlineskip
\vspace*{12pt}
With some extra notation, we can state the theorem more succinctly.  Let
${\cal V}_i$ be the vector space of Vassiliev invariants spanned by $v_i$ from
$(\ref{exp})$.  Consider the underlying vector space of
the polynomial algebra ${\bf Q}[v]$ with basis $\{v^0,v^1,v^2,\ldots\}$.  Let
$E$ be the vector space isomorphism $E: {\bf Q}[v] \to \oplus {\cal V}_i$, 
where $E(v^i) = i!v_i$.  We therefore obtain:
\[ a_n = E(f_{d,n}(v)) \]

\begin{example} \rm As in Section 2, let $\alpha_{d,n}$ denote the $n^{th}$
coefficient of the Jones polynomial of degree $\leq d$.
\[ \begin{array}{rcl}
\alpha_{2,0} &=& \frac{1}{2!2!} (4-5v_2+v_4) = \frac{1}{2!2!} E \left(
(v-2)(v-1)(v+1)(v+2) \right) \\ \\
\alpha_{3,0} &=& \frac{1}{3!3!} (36-49v_2+14v_4-v_6) = \\ \\ 
         &=& -\frac{1}{3!3!} E \left( (v-3)(v-2)(v-1)(v+1)(v+2)(v+3)
\right) \\ \\
\alpha_{2,1} &=& \frac{1}{1!3!} (4v_1+4v_2-v_3-v_4) = -\frac{1}{1!3!} E \left(
(v-2)(v)(v+1)(v+2) \right) \\ \\
\alpha_{3,1} &=& \frac{1}{2!4!} (36v_1+36v_2-13v_3-13v_4+v_5+v_6) = \\ \\ 
         &=& \frac{1}{2!4!} E \left( (v-3)(v-2)(v)(v+1)(v+2)(v+3) \right) \\ \\ 
\end{array} \]
\end{example}

\begin{remark} \rm
From Theorem \ref{Jb}, we obtain a formula for any $v_i(L)$ in terms of $v_0(L),\dots,v_{2d}(L)$.
Namely, 
\[ v_i(L)= \frac{1}{i!} \sum_{k=-m}^{n} k^ia_k(L)
=\sum_{j=0}^{2d} \left(\frac{1}{i!} \sum_{k=-m}^{n} k^if_{d,k}^{(j)}(0)\right) v_j(L) \]
Similarly, in \cite{KSS} it was shown, without an explicit formula, that for any knot $L$, 
there exists $N,$ such that all $v_i(L)$ are determined by $v_0(L),\dots,v_N(L)$.

\end{remark}

\section{Approximating Jones coefficients by Vassiliev invariants:
Infinite case}

In this section, we formally let $d \to \infty$ to find the correct
formula for the coefficients of an arbitrary Jones polynomial of a knot, and then
prove that the resulting series of Vassiliev invariants converges.
We also extend the vector space isomorphism 
$E: {\bf Q}[[v]] \to \prod {\cal V}_i$, where $E(v^i) = i!v_i$.

We first consider $a_0$:
\[ f_{d,0}(v) = \prod_{{j \neq 0}\atop {j=-d}}^{d} \frac{v-j}{-j} =
\prod_{{j \neq 0}\atop {j=-d}}^{d} \left( 1 - \frac{v}{j} \right) =
\prod_{j=1}^{d} \left( 1 - \frac{v}{j} \right) \left( 1 + \frac{v}{j}
\right) = \prod_{j=1}^{d} \left( 1 - \frac{v^2}{j^2} \right) \]

We recall the Weierstrass product factorization for entire functions:
\[ \prod_{j=1}^{\infty} \left( 1 - \frac{z^2}{j^2} \right) = \frac{\sin
\pi z}{\pi z} \]

We define $f_{\infty,0}(v) = \lim_{d \to \infty} f_{d,0}(v)$, so we obtain
\[f_{\infty,0}(v) = \frac{\sin \pi v}{\pi v} = 1 - \frac{{(\pi v)}^2}{3!}
+ \frac{{(\pi v)}^4}{5!} -\frac{{(\pi v)}^6}{7!} + \cdots \]

Formally (we prove convergence below), we obtain the following beautiful
formula:
\begin{equation}
 a_0 = E(f_{\infty,0}(v)) = E(\frac{\sin \pi v}{\pi v}) =  v_0 -
\frac{\pi^2}{3}v_2 + \frac{\pi^4}{5}v_4 -\frac{\pi^6}{7}v_6 + \cdots
\end{equation}

We now consider $a_n$:
\[ f_{d,n}(v) = \prod_{j\neq n \atop j=-d}^{d}  \frac{v-j}{n-j} =
\prod_{j=-d}^{n-1}\frac{v-j}{n-j}\prod_{j=n+1}^{d}\frac{v-j}{n-j} \]
Let $k=n-j$ in the first product, and $k=j-n$ in the second product, so we
obtain
\[ f_{d,n}(v) = \prod_{k=1}^{d+n}\frac{v+k-n}{k}
\prod_{k=1}^{d-n}\frac{v-k-n}{-k}=
\prod_{k=1}^{d-n}(1 - \frac{(v-n)^2}{k^2}) \prod_{k=d-n+1}^{d+n}(1 +
\frac{v-n}{k}) \]

The second product is finite:  let $l=k-d$, then we obtain
$\prod_{l=-n+1}^{n}(1 + \frac{v-n}{l+d})$.  For any $n$, as $d \to \infty$
we can easily see that this product converges to $1$.  The first product
converges for all $n$ by the same argument as above.  This suggests the following theorem:
\textlineskip
\vspace*{12pt}

\begin{thm} \label{T}
\[ a_n = E(f_{\infty,n}(v)) = \sum_{i=0}^{\infty}f_{\infty,n}^{(i)}(0)v_i,\;
{\rm  where }\; f_{\infty,n}(v) = \cases{ \frac{\sin \pi (v-n)}{\pi(v-n)} & if $v \neq n$ \cr
 \qquad 1 & if  $v = n$ } \]
\end{thm}

{\bf Proof.}
For any given knot $L$, the Jones polynomial has finite degree $d$.
\[ \sum_{i=0}^{\infty}f_{\infty,n}^{(i)}(0)v_i(L) =
\sum_{i=0}^{\infty}f_{\infty,n}^{(i)}(0)\left( \frac{1}{i!}
\sum_{k=-d}^{d}k^ia_k(L) \right) \]
\[ =
\sum_{k=-d}^{d}a_k(L)\left(\sum_{i=0}^{\infty}\frac{1}{i!}f_{\infty,n}^{(i)}(0)k^i\right)\]
\[ = \sum_{k=-d}^{d}a_k(L)f_{\infty,n}(k) = a_n(L) \]
\hfill \qed\kern0.8pt

\begin{example} 
\[ a_1  =  v_1 + 2v_2 + (3!-\pi^2)v_3 + (4!-4\pi^2)v_4
+(5!+\pi^4-20\pi^2)v_5 + \cdots \]
\end{example}

\begin{corollary}
Any Jones coefficient of a knot can be approximated by Vassiliev invariants $\tilde{v}_i$ of order $i$:
\[ a_n = \lim_{i \to \infty} \tilde{v}_i,\quad  {\rm where }\;
\tilde{v}_i = \sum_{j=0}^{i}f_{\infty,n}^{(j)}(0)v_j \]
\end{corollary}

\begin{remark} \rm
Jones coefficients can be shown not to be Vassiliev invariants by considering twist sequences \cite{T}.  Let $T_m$ be the $(2,2m+1)$-torus knot. By Theorem 2.2.1 of \cite{T}, the restriction of any Vassiliev invariant to the sequence $\{T_m\}$ is a polynomial in $m$.  However, $J_{T_m}(t)= -t^m(t^{2m+1}-\ldots+t^3-t^2-1)$.  Thus for $n \geq 0$, $a_n(T_m)=0$ for $m<\frac{n-1}{3}$ or $m>n$.  If $a_n(T_m)$ were a polynomial in $m$, it would be zero on $T_m$, which is clearly false.  Similarly, we can take mirror images of $T_m$ to show that for all $n<0$, $a_n$ is not a Vassiliev invariant.  (See also \cite{Z}.)

Trapp \cite{T} also showed that if a sequence of Vassiliev invariants converges
uniformly for all knots, then the limit is also of finite type.  Because
the Jones coefficients are not of finite type, the pointwise
limits above cannot be uniformly convergent for all knots.  Indeed, the proof of
Theorem $\ref{T}$ requires us to first choose a particular knot.
\end{remark}
\textlineskip
\vspace*{12pt}

\begin{remark} \rm
The function $f_{\infty,n}(v)$ is not unique, because the infinite
Vandermonde matrix can have infinitely many left inverses.  Since every
Jones polynomial has finite degree $d$, but Vassiliev invariants may be
nonzero for arbitrary orders, we can view the infinite matrix as a linear
operator $\oplus_{i=-\infty}^{\infty} {\bf Q} \to \prod_{i=0}^{\infty} {\bf Q}$, so
it can have infinitely many left inverses, but no right inverse.
The proof of Theorem $\ref{T}$ only requires that $f_{\infty,n}(v)$ has a
Taylor expansion about zero, and that $f_{\infty,n}(m)=\delta_{m,n}$, the
Kronecker pairing for all $m, n \in {\bf Z}$.  If we also insist that
$f_{\infty,n}(v) = \lim_{d \to \infty} f_{d,n}(v)$, then the generating
function depends on how we select invertible finite blocks to exhaust the
infinite Vandermonde matrix.  
\end{remark}
\textlineskip
\vspace*{12pt}

\begin{remark} \rm
Given an infinite sequence $\{v_n(L)\}$, we can also approximate the
degree $d$ of the Jones polynomial by functions of finite type invariants:
\[ d = \lim_{n \to \infty} \sqrt[n]{\vline \sum_{k=-d}^{d} k^n
a_k(L) \vline } = \lim_{n \to \infty} \sqrt[n]{n!|v_n|}\]
\end{remark}

\begin{remark} \rm
All of the results above carry over with slight modifications to links.
The Jones polynomial of a link may be a Laurent polynomial times $\sqrt{t}$, so
instead of (\ref{Van}), the matrix and resulting formulas appear as a special 
case of Theorem \ref{H}.
\end{remark}

\section{Approximations of coefficients of specializations of HOMFLY and
Kauffman polynomials}

The Jones polynomial is a specialization of both the HOMFLY and Kauffman
two-variable polynomials, $H_L(a,z)$ and $F_L(a,z) \in
{\bf Z}[a^{\pm1},z^{\pm1}]$.  We consider an infinite sequence of one-variable
specializations of the HOMFLY polynomial.  The same result and proof
applies to the Kauffman polynomial as well.  Let $N \in
{\bf Z}\backslash\{0\}$.
\begin{equation} \label{Haz}
 a=t^{N/2}, z=t^{1/2}-t^{-1/2} \quad \Rightarrow \quad H^N_L(t) \in
Z[t^{\pm\frac{1}{2}}]
\end{equation}

\noindent
The Jones polynomial is obtained at $N=-2$.  Now, suppose 
$ H^N_L(t)=b^N_{-m}t^{-m/2} + \cdots + b^N_0 + \cdots + b^N_nt^{n/2}. $
Let $d=\max(m,n)$.  In \cite{BL}, it was shown that if we let $t=e^x$,
\begin{equation}
H^N_L(e^x)  = \sum_{i=0}^{\infty} \left( \frac{1}{i!}
\sum_{k=-m}^{n}{\left(\frac{k}{2}\right)}^ib^N_k(L) \right)x^i =
\sum_{i=0}^{\infty} v^N_i(L)x^i
\end{equation}
then $v^N_i$ is a Vassiliev invariant of order $i$.
\textlineskip
\vspace*{12pt}

\begin{thm} \label{H}
For any link $L$ and any $N \in {\bf Z}\backslash\{0\}$, let its $N^{th}$
HOMFLY polynomial be $H^N_L(t) = \sum_{n=-d}^d b^N_n t^{n/2}$.  Let
$E^N(v^i)= i!v^N_i$.  Then, 
\[ b^N_n = E^N(f_{d,n}(v)) = \sum_{i=0}^{2d} f_{d,n}^{(i)}(0) v^N_i, \quad
{\rm  where }\; f_{d,n}(v) = \prod_{j\neq n \atop j=-d}^{d} 
\frac{2v-j}{n-j} \]
\[ b^N_n = E^N(f_{\infty,n}(v)) =
\sum_{i=0}^{\infty}f_{\infty,n}^{(i)}(0)v^N_i,\; {\rm where }\;
f_{\infty,n}(v) = \cases{ \frac{\sin \pi (2v-n)}{\pi (2v-n)} & if  $2v \neq n$ \cr
 \qquad 1 & if  $2v = n$ } \]
\end{thm}

{\bf Proof.}
The proof is just a modification of the proof for the Jones polynomial.
Instead of (\ref{Van}), we have
\[
  	\left( \matrix{
	\cdots&1&1&1&1&1&\cdots \cr  
	\cdots&-1&-1/2&0&1/2&1&\cdots \cr 
	\cdots&(-1)^2&(-1/2)^2&0&(1/2)^2&1^2&\cdots \cr 
	&\vdots&\vdots&\vdots&\vdots&\vdots \cr 
	} \right) 
	\left( \matrix{
	\vdots \cr b^N_{-1} \cr b^N_0 \cr b^N_1 \cr  \vdots \cr 
	} \right)
=
	\left( \matrix{
	v^N_0 \cr v^N_1 \cr 2!v^N_2 \cr 3!v^N_3 \cr  \vdots \cr 
	} \right)	
\]
To find the inverse of both the infinite matrix and any invertible finite
block, we just need the following:
\textlineskip
\vspace*{12pt}

\begin{proposition} 
For any $m \in {\bf Z}$ such that $-d \leq m \leq d$,
\[ f_{d,n} ( \frac{m}{2} ) = \cases{ 1 & if  $m=n$ \cr 0 & if  $m \neq n$ } \]
and similarly for $f_{\infty,n}(\frac{m}{2})$ for any $m \in {\bf Z}$. \hfill \qed\kern0.8pt
\end{proposition}
\textlineskip
\vspace*{12pt}

\begin{corollary} \label{Hz}
For any link $L$, any fixed complex number $z$, and any $N \in
{\bf Z}\backslash\{0\}$, $H^N_L(z)$ is a limit of Vassiliev invariants.
\end{corollary}
\textlineskip
\vspace*{12pt}

{\bf Proof.}
The proof is the same as the proof of Corollary \ref{z}.
\hfill \qed\kern0.8pt
\textlineskip
\vspace*{12pt}

\begin{remark} \label{Faz} \rm
Following \cite{BL}, let $F_L(a,z)$ denote the Dubrovnik version of the
Kauffman polynomial.  To obtain one-variable specializations, for $N \in
{\bf Z}\backslash\{0\}$, set $a=t^N, z=t-t^{-1}$ and let $F^N_L(t)$ denote the
$N^{th}$ Kauffman polynomial.  Theorem \ref{H} and Corollary \ref{Hz}
carry over to $F^N_L(t)$.
\end{remark}

\section{Approximations of other link invariants}

Let $\Sigma_n(L)$ denote the n-fold branched cover of a link $L$.  
Let $Q_L(x)$ be the specialization at $a=1$ of the standard Kauffman
polynomial.  
This link polynomial satisfies the skein relation
$Q_{L_+}+Q_{L_-}=x(Q_{L_0} + Q_{L_{\infty}})$
and $Q(\mbox{unknot})=1$ \cite{BLM, Ho}.
A lot of information about $H_1(\Sigma_n(L))$ can be obtained by
evaluating
link polynomials of $L$ at special values. In fact, it seems reasonable to 
conjecture that $H_1(\Sigma_2(L), {\bf Z})$ is determined by $Q_L$ \cite{Rong}.
We summarize the results below.  All the evaluations can be found in
\cite{Lickorish}, except for $Q_L(2\cos \frac{2\pi}{5})$ which is
given in \cite{Jones} and \cite{Rong}.

For the Jones polynomial, the most interesting values 
to evaluate are $t=e^{\frac{2\pi i}{r}}$, where $r$ is a positive
integer.

\begin{table}[htbp]
\tcaption{Jones polynomial at $t=e^{\frac{2\pi i}{r}}$}
\centerline{\footnotesize\smalllineskip
\begin{tabular}{c|l}\\
$r$ & $J_L(e^{\frac{2\pi i}{r}})$  \\ \hline
$r=1$ & $(-2)^{\ell-1} $ \\ \hline
$r=2$ & $|H_1(\Sigma_2, {\bf Z})|=\det L$ \\ \hline
$r=3$ & $1$ \\ \hline
$r=4$ & $(-\sqrt{2})^{(\ell-1)}(-1)^{\mbox{Arf}(L)} $  if Arf$(L)$
exists,\\
 & $0$ if Arf$(L)$ undefined \\ \hline
$r=6$ & 
$(-\sqrt{3})^{\dim (H_1(\Sigma_2, {\bf Z}_3))}(-i)^w $ 
\\
\end{tabular}}
\end{table}

\noindent
In the table, $\ell$ is the number of components of $L$,  
Arf$(L)$ is the Arf invariant of $L$, and $w\in {\bf Z}_4$ is the Witt class of
the Seifert form mod 3 of $L$.
 
Because the Jones polynomial is a specialization of the HOMFLY polynomial, 
the table also gives evaluations of the HOMFLY polynomial.
Another interesting value not listed above is 
$H_L(-i,i)=(i\sqrt{2})^{\dim (H_1(\Sigma_3, {\bf Z}_2))}$.
Let us also recall that 
$|H_1(\Sigma_n, {\bf Z})| = |\prod \Delta_L(r_i)|$,  where $r_i$'s are the 
$n$th roots of unity.

For the $Q$-polynomial, the interesting values are at $x=2\cos
\frac{2\pi}{r}=q+q^{-1}$, where $q=e^{\frac{2\pi i}{r}}$.

\begin{table}[htbp]
\tcaption{$Q$-polynomial at $x=2\cos\frac{2\pi}{r}$}
\centerline{\footnotesize\smalllineskip
\begin{tabular}{c|l}\\
 $r$ & $Q_L(2\cos \frac{2\pi}{r}) $ \\ \hline
$r=1$ & $(-1)^{\ell-1}|\det L|^2 $ \\ \hline
$r=2$ & $(-2)^{\ell-1}$ \\ \hline
$r=3$ & $(-3)^{\dim (H_1(\Sigma_2, {\bf Z}_3))}$ \\ \hline
$r=4$ & undefined \\ \hline
$r=5$ & $(-1)^{t_5(L)}(\sqrt{5})^{\dim (H_1(\Sigma_2, {\bf Z}_5))}$  \\ \hline
$r=6$ & 1 \\
\end{tabular}}
\end{table}

\noindent
Here $t_5$ is $0$ or $1$,  and can be written in terms of the Seifert form
mod $5$ \cite{Rong}.  
\textlineskip
\vspace*{12pt}

\begin{thm}
Let $\phi$ be any of the following knot invariants:
$\det L,  |H_1(\Sigma_2, {\bf Z}_p)|$ where $p=3, 5$, $|H_1(\Sigma_3, {\bf Z}_2)|$,
$|H_1(\Sigma_n,{\bf Z})|$.  Then\\
(a) $\phi$ is not a Vassiliev invariant.\\
(b) $\phi$ is a limit of functions of Vassiliev invariants.
\end{thm}


{\bf Proof.}
(a) As remarked in \cite{ST}, a knot invariant $\phi$ is not
a Vassiliev invariant if there is a knot $K$ with $\phi(K\# K')\neq 
\phi(\mbox{unknot})$ for all knots $K'$.
Now let $\phi$ be, for example, the order of $H_1(\Sigma_2, {\bf Z}_p)$.
Let $K$ be a knot with $\phi(K)=p$ (e.g., any 2-bridge knot for which 
$\Sigma_2(K)$ is the lens space $L_{p,q}$).
Then $\phi(K\#K')=\phi(K)\phi(K')\neq \phi(\mbox{unknot})$.
Similarly, $|H_1(\Sigma_3, {\bf Z}_2)|$ and $|H_1(\Sigma_n,{\bf Z})|$ are not
Vassiliev invariants.

(b) The following equations come from the tables and comments above. 
Together with Corollary \ref{z} and Corollary \ref{Hz}, they imply
that $\phi$ is a limit of functions of Vassiliev invariants. 

\[ \begin{array}{lcl}
|H_1(\Sigma_2, {\bf Z})| &=& J_L(-1) \\
|H_1(\Sigma_2, {\bf Z}_3)| &=& |J_L(e^{\frac{\pi i}{3}})|^2 \\
|H_1(\Sigma_2, {\bf Z}_5)| &=&  |Q_L(2 \cos(2\pi/5))|^2 \\
|H_1(\Sigma_3, {\bf Z}_2)| &=& |H_L(-i,i)|^2 \\
|H_1(\Sigma_n, {\bf Z})| &=& |\prod \Delta_L(r_i)|, \mbox{ where 
      $r_i$'s are the $n$th roots of unity.}
\end{array} \]

Note from (\ref{Haz}), we obtain that $H_L(-i,i)$ is $H^N_L(t)$, where
$N=9, \; t=e^{i\pi/3}$.

As in Section 5, let $F^N_L(t)$ denote the $N^{th}$ Dubrovnik Kauffman
polynomial.  By a change of variables, $Q_L(x)=(-1)^\ell F_L(i,-ix)$
\cite{K}.  Thus, by Remark \ref{Faz} we obtain that $Q_L(2\cos(2\pi/5))$ is
$(-1)^\ell F^N_L(t)$, where $N= -5, \; t=-ie^{2i\pi/5}$.

Since the coefficients of the Alexander-Conway polynomial are Vassiliev
invariants \cite{B}, by an argument similar to Corollary \ref{z} we obtain 
that $\prod \Delta_L(r_i)$ is a limit of Vassiliev invariants.
Thus, $|\prod \Delta_L(r_i)|$ is a limit of the absolute value function
 of Vassiliev invariants.
\hfill \qed\kern0.8pt
\textlineskip
\vspace*{12pt}

\begin{remark} \rm
For all $\phi$ except for $\phi=|H_1(\Sigma_n, {\bf Z})|$, we can show that
$\phi$ is actually a limit of Vassiliev invariants.  For example,
in the case of $|H_1(\Sigma_2, {\bf Z}_3)|$, it follows from Corollary \ref{z} 
and the equation: 
$|H_1(\Sigma_2, {\bf Z}_3)| = J_L(e^{\frac{\pi i}{3}}) J_L(e^{-\frac{\pi i}{3}})$,
since the two factors on the right are complex conjugates.

For other functions, e.g., $\phi=\dim (H_1(\Sigma_2, {\bf Z}_3))$, 
the argument above can only show that $\phi$ is a limit of
functions of Vassiliev invariants.  Note that in the case of $|H_1(\Sigma_n, {\bf Z})|$, 
the $n$-th roots of unity all appear in conjugate pairs, thus the product 
$\prod \Delta_L(r_i)$ can be arranged in conjugate pairs. It follows that 
the absolute value sign is not needed, except for $\Delta(-1)\Delta(1)$.
\end{remark}
\textlineskip
\vspace*{12pt}

\begin{remark} \rm
Another knot invariant which is not a Vassiliev invariant, but is a limit of Vassiliev invariants is $tri(L)$, the number of 3-colorings of $L$.  This follows from a result of Przytycki \cite{P1, P2}: $tri(L) = 3|J_L(e^{\frac{\pi i}{3}})|^2$.
\end{remark}

\section{Conclusion}

Here, we make some final remarks on approximations by Vassiliev
invariants.
Our work is motivated by the following two equivalent conjectures (see \cite{Rog}):
\pagebreak

\noindent
{\bf Conjecture 7.1.} Vassiliev invariants separate knots.  That is,
for any two knots $K_1$ and $K_2$, there is a Vassiliev invariant
$v$ with $v(K_1)\neq v(K_2)$.
\textlineskip
\vspace*{12pt}

\noindent
{\bf Conjecture 7.2.}  Every knot invariant is a limit of Vassiliev
invariants.
That is, for any knot invariant $\phi$, there is a sequence of Vassiliev 
invariants $v_n$ with $\phi(L)=\lim_{n\rightarrow \infty} v_n (L)$ for all
knots $L$.
\textlineskip
\vspace*{12pt}

We have verified that a number of knot invariants (e.g., coefficients of
link polynomials) are indeed limits of Vassiliev invariants.  Some other
knot invariants (e.g., the degree of the Jones polynomial) are proved to
be limits of functions of Vassiliev invariants. In light of this, we
propose:
\textlineskip
\vspace*{12pt}

\noindent
{\bf Conjecture 7.3.}  Every knot invariant is a limit of functions of
Vassiliev invariants.  That is, for any knot invariant $\phi$, there is a 
sequence of functions $f_n$ and Vassiliev invariants $v_n$ with 
$\phi(L)=\lim_{n\rightarrow \infty} f_n(v_n) (L)$ for all knots $L$.
\textlineskip
\vspace*{12pt}

In general, if $f$ is an analytic function, $v$ is a Vassiliev invariant,
then it is not hard to show that $f(v)$ is a limit of Vassiliev
invariants.
Consequently, a limit of analytic functions of Vassiliev invariants is in
fact
a limit of Vassiliev invariants. 
However, this is not clear if the functions are not analytic functions.
This is the case, for example, for the degree of the Jones polynomial,
where the functions $f_n$ are the $n$th root function.

One good aspect of Conjecture 7.3 is that it is easier to verify than
Conjecture 7.2 for a given knot invariant, but is still strong enough to imply Conjecture 7.1.
Therefore by \cite{Rog}, 
\[ {\rm Conjecture \,7.2} \Rightarrow {\rm Conjecture \,7.3} \Rightarrow
{\rm Conjecture \,7.1} \Rightarrow {\rm Conjecture \,7.2}. \]

\nonumsection{Acknowledgements}
We would like to thank Dror Bar-Natan, Xiao-Song Lin, Ted Stanford, and Ed Swartz for helpful discussions.  The first author was partially supported by NSF grant DMS-98-03518.  The second author was partially supported by NSF grant DMS-97-29992 while visiting the Institute for Advanced Study.

\nonumsection{References}

\end{document}